\documentstyle[amsfonts,12pt,epsfrag,psfrag]{amsart}

\newtheorem{lemma}{Lemma}[section]

\newtheorem{theorem}{Theorem}
\newtheorem{conjecture}{Conjecture}
\newtheorem{corollary}[lemma]{Corollary}

\newcommand{\bbR}{\mbox{${\bf I\! R}$}}
\newcommand{\Rt}{\mbox{${\bbR}^3$}}

\begin{document}

\title{Vertex Theorems for Capillary Drops on Support Planes}
\author{Robert Finn \and John McCuan}
\address{\hskip-\parindent
Robert Finn, Mathematics Department, Stanford University, Stanford, 
CA  94305-2125} 
\email{finn@@math.stanford.edu}
\address{\hskip-\parindent
John McCuan, Mathematical Sciences Research Institute, 1000 
Centennial Drive, Berkeley, CA  94720-5070}
\email{john@@msri.org}
\thanks{\tiny This work was supported in part by a grant from the 
National Aeronautics and Space Administration, and in part by a grant from the 
National Science Foundation. Part of the work was completed while the 
latter author was a National Science Foundation Postdoctoral Fellow at the 
Mathematical Sciences Research Institute in Berkeley, CA.  The authors 
wish to thank Universit\"at Leipzig for its hospitality during the 
initial phases of this work, and the Max-Planck Institute f\"ur 
Mathematik in Leipzig for its hospitality during completion of the 
work}

\begin{abstract}
We consider a capillary drop that contacts several planar bounding walls 
so as to produce singularities (vertices) in the boundary of its free 
surface.  It is shown under various conditions that when the number of 
vertices is less than or equal to 
three, then the free surface must be a portion of a sphere.  These results 
extend the classical theorem of H. Hopf on constant mean curvature 
immersions of the sphere.  The 
conclusion of sphericity cannot be extended to more than three 
vertices, as we show by examples.
\end{abstract}



\maketitle

\section{Overview}

We consider in this paper liquid drops in  $\Rt$  resting in mechanical 
equilibrium in the absence of gravity on rigid support surfaces that 
consist of a finite number  of intersecting planes $\Pi_j$.  We restrict 
attention 
to the physically familiar configurations for which the surface interface 
${\cal S}$ is topologically a disk.  Typical configurations with 
four and eight vertices are indicated in Figure 1. 
\begin{figure}[hb]
\centerline{{\psfrag{x}{$x$}
                   \psfrag{S}{${\cal S}$}
                   \psfrag{t}{$\gamma$}
                   \psfrag{p1}{$\Pi_1$}
                   \psfrag{p2}{$\Pi_2$}
         \epsfysize=3.6in
         \epsfxsize=3.6in
         \leavevmode\epsfbox{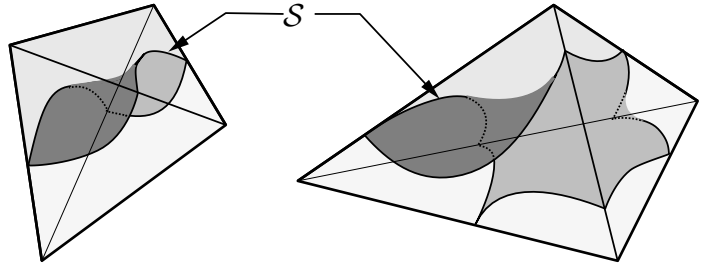}}}
\centerline{ }
\caption{Typical Configurations\label{typical}}
\end{figure}

Under the indicated conditions, ${\cal S}$ will have constant mean 
curvature $H$, 
and it is natural to ask to what extent a theorem of H.~Hopf \cite{HopUbe} 
on closed immersions of constant $H$ will apply to the configurations 
we consider. That is, we wish to determine criteria under which ${\cal S}$ 
will necessarily be part of a metric sphere. We will present such 
conditions, and we will provide examples to show that our criteria are 
reasonably sharp. 

Our results relate closely to (and in fact were inspired by) recent 
discoveries on tubular liquid bridges \cite{McCSym} and on 
capillary surfaces 
in cylindrical tubes with protruding edges \cite{ConCap, FinLoc}; 
these papers form 
the background for the perspective that we adopt. 
The support surfaces ${\cal P}$ considered here consist of a finite 
collection of 
planes, no three of which intersect in a common line, and which together 
contain the boundary of an open connected region ${\cal I}$ in $\Rt$.  
We assume ${\cal S}$ to be a constant mean curvature surface, lying 
in the closure of ${\cal I}$, whose trace on each supporting plane 
$\Pi_j$ is a smooth {\em contact line} ${\cal C}_j$ along which ${\cal S}$ 
meets $\Pi_j$ in a constant {\em contact angle} $\gamma_j$.  Formally, 
${\cal S}$ satisfies the variational condition
$$
	\delta\left(|{\cal S}| + \sum\beta_j{\cal S}_j +2HV\right) = 0, 
$$ 
with $\beta_j = \cos\gamma_j$, ${\cal S}_j$ the wetted area on the plane 
$\Pi_j$, and $V$ the drop volume, see, e.g., \cite[Chapter 1]{FinEqu}. 
However, much of our formal analysis will encompass immersions ${\cal S}$ 
for which the enclosed volume may not be defined, and we focus attention 
on solutions of the system of (Euler-Lagrange) equations
$$
	\Delta\vec{x} = 2H\vec{N}
$$
under the indicated boundary conditions.  Here, the Laplacian is the 
intrinsic operator on the surface, and $\vec{N}$ is a unit normal to 
${\cal S}$.

Physically, the constancy of the $\gamma_j$ 
means that the liquid and each of the planes are assumed to be of 
homogeneous materials. The materials (and hence the angles $\gamma_j$) 
may differ 
from plane to plane, but the condition is introduced that the angle pair 
$(\gamma_j,\gamma_k)$ on any two adjacent planes in contact with the drop, 
and for which the intersection line ${\cal L}_{jk}$ passes through 
${\cal S}$, lies interior to a certain rectangle ${\cal Q}$, see 
Section 2 below. We assume also that the angles between any two 
intersecting support planes, measured in an appropriate sense, 
are less than 
$\pi$. There is strong heuristic evidence that both these conditions 
are (in general) necessary. No boundary condition is introduced on the 
lines ${\cal L}_{jk}$ themselves.  With regard to behavior near these 
lines, our 
conceptual point of view will be that any contact point of ${\cal S}$ 
with an ${\cal L}_{jk}$ is a {\em vertex} ${\cal V}$.  Vertices will turn out to 
be uniquely determined points, 
but it is possible (and we think it desirable) to start off in some cases 
with weaker hypotheses, which do not initially require ${\cal S}$  
to be defined on the ${\cal L}_{jk}$.  

A case of primary interest is the dihedral angle consisting of 
two planes $\Pi_1$ and $\Pi_2$ that intersect in a line ${\cal L}$.  
We discuss this case heuristically in Section~2 below.  We enumerate 
in Hypothesis~A physically natural conditions 
applying to a drop in a dihedral angle that serve as 
motivation for the weaker conditions of Hypothesis~B below.

\medskip

\noindent {\bf Hypothesis A.}   {\em ${\cal S}$ is  globally embedded, and 
together with portions of the supporting planes bounds a drop volume 
topologically a ball.  Each 
vertex ${\cal V}$ is a unique point on ${\cal L}$, 
and ${\cal S}$ can 
be parametrized locally up to ${\cal V}$ by  continuous  functions. 
There exists a plane  $\Pi$ orthogonal to ${\cal L}$, 
cutting off a 
portion ${\cal S}_{\cal V}$ of ${\cal S}$ containing ${\cal V}$, such 
that the reflection of ${\cal S}_{\cal V}$ in that plane lies interior 
to the drop volume. }

\medskip

We will now state precisely several more general assumptions that are 
sufficient to imply our main results and are used (sometimes without 
explicit mention) throughout the paper.  We have not attempted to find 
the weakest possible assumptions except for those directly related to 
the main contribution of the paper, vis. the singular behavior of the 
boundary curve.  The following four conditions collectively 
generalize Hypothesis~A. 

\medskip

\noindent {\bf Hypothesis~B.}

\smallskip

Let $D$ be the closed unit disk in $\bbR^2$ and $v_i$, 
$i=1,\ldots, V$ a finite clockwise-ordered collection of points in 
$\partial{D}$.  For $i=1,\ldots, V$, let $A_i$ be the arc on $\partial{D}$ 
between $v_i$ and $v_{i+1}$; $A=\cup A_i$.  

\smallskip

\noindent {\small\bf 1. Topological Condition:}  {\em There is a local 
homeomorphism $\phi$ 
of $D_v = D\backslash \{v_i\}$ onto ${\cal S}$, i.e., for each $a\in D_v$ 
there is some neighborhood $B(a)$ of $a$ in $D_v$ such that $\phi$ restricted 
to $B(a)$ is a homeomorphism.}

\smallskip

The following subsets of ${\cal P}$ have special significance:  If $\Pi_j$ 
and $\Pi_k$ intersect in a line ${\cal L}_{jk}$, then ${\cal E}_{jk} =
{\cal L}_{jk}\cap \bar{{\cal I}}$ is the {\em $jk$-edge.}\/  
The {\em edge of ${\cal P}$}\/ is ${\cal E} = \cup {\cal E}_{jk}$.  
The {\em faces} of ${\cal P}$ are ${\cal F} = \bar{{\cal I}} \backslash 
{\cal E}$.  No confusion results if we refer to a 
connected component of ${\cal F}$ as $\Pi_j$ or $\Pi_k$.  

\smallskip

\noindent{\small\bf 2. Smoothness Condition:}  {\em $\phi_{|_A}$ can be made 
locally smooth, i.e., for each $a\in A$ there is a homeomorphism 
$\psi: (-1,1) \to A$ such that $\psi(0) = a$ and $X^{\partial} =
\phi\circ\psi: (-1,1) \to {\cal F}$ represents a differentiable 
embedded curve.  

Also, $\phi_{|_{D_v}    }$ can be made locally smooth, i.e., for each $a \in 
{\rm int}\, D$ there is a homeomorphism 
$\psi: B_1(0)\to {\rm int}\, D$ 
such that $X=\phi\circ\psi: B_1(0)\to {\cal I}$ represents a $C^2$ (open) 
embedded surface, and for $a\in A$ a similar statement holds with $\psi$ 
defined on a half neighborhood 
$B^+_1(0) = \{(a_1,a_2)\in B_1(0): a_2\ge 0\}$.}

\smallskip

\noindent Notice that we have required each boundary component 
${\cal C}_i = \phi(A_i)$ to lie in some face $\Pi_j$;  it is not required 
that this association be one-to-one.

\smallskip

For each pair of planes $\Pi_j$, $\Pi_k$ in ${\cal P}$ that intersect in 
a line ${\cal L}_{jk}$, denote by $\Pi =  \Pi_{jk}$ a plane orthogonal 
to $\Pi_j$ and $\Pi_k$.  The heart of Hypothesis B is the following

\smallskip

\noindent {\small\bf 3. Vertex Condition:}   {\em   To each $v_i\in \partial D$, 
there is associated a pair of intersecting planes 
$\Pi_j$, $\Pi_k\in {\cal P}$ and two neighborhoods: $B$, a neighborhood of 
$v_i$ in $D$, and ${\cal N}$, a neighborhood of $O= {\cal L}_{ij} \cap \Pi$ in 
$\Pi \cap \bar{{\cal I}}$, such that $\phi(B\backslash v_i)$ is a graph 
$u(p)$ over ${\cal N} \backslash O$.}

\medskip

\noindent In the case of Hypothesis B, each such triple 
$\{{\rm plane pair},\ {\rm neighborhood},\ {\rm function}\} = 
\{ (\Pi_j,\Pi_k), {\cal N}, u \}$ will be 
said to determine a vertex  ${\cal V}$, and we write 
$S_{\cal V} = \phi(B\backslash v_i)$, ${\cal N}_{\cal V}= {\cal N}$, 
$u_{\scriptscriptstyle {\cal V}}= u$, etc..  As noted above, vertices 
need not in this case be initially defined as points;  we shall 
however prove the existence of the vertex points and the smoothness of
$u_{\scriptscriptstyle\cal V}(p)$  up to those points, and in fact we 
shall do so without 
growth hypotheses on the graph.  Once that is done, the local 
homeomorphism $\phi$ that parameterizes ${\cal S}$ can 
be extended continuously to $D$.   

\smallskip

In general however, we require an additional condition in order to 
separate the vertices one from another.  
Though not the weakest possible condition, in order to simplify the 
proofs of Theorem~2 and of Section~5 below we impose the following 

\smallskip

\noindent {\small\bf 4. Separation Condition:} {\em For each vertex 
${\cal V}$, 
${\cal S}_{\cal V}$ can be chosen so that 
$\phi: \phi^{-1}({\cal S}_{\cal V}) \to {\cal S}_{\cal V}$ is a local 
homeomorphism, and if there is another vertex ${\cal V}'$ corresponding 
to the same plane pair, then there exists a sequence 
$p_j\to {\cal N}_{\cal V}\cap {\cal L}_{jk}$ on $\Pi$, on which 
$|u_{\scriptscriptstyle\cal V}-u_{\scriptscriptstyle\cal V'}|$ 
is bounded away from zero.}


\smallskip

\noindent If $\phi$ extends continuously to the boundary, then the 
Separation Condition 
becomes simply:  {\em For each vertex ${\cal V}$, 
$\phi^{-1}({\cal V})= v_i$} 



\medskip

Under 
either of the above hypotheses, and under the 
above conditions on the opening angles and the adjacent contact angle 
pairs, we intend to prove that if $V$  is the number of vertices and if 
$V\le 2$ then ${\cal S}$ is a portion of a metric sphere. We will 
obtain the same result when $V = 3$, under an additional condition on 
the orientation of ${\cal S}$.  Remarkably, the orientation required is 
the reverse of the (uniquely determined) 
orientation that occurs when $V \le 2$.  The interest in these results is 
underscored by the fact that if $V\ge 4$ then ${\cal S}$ need not be 
spherical, as we show by example.

	If we restrict attention to a dihedral angle and the case 
$V = 2$ (drop in a wedge), 
then the Vertex and Separation Conditions can be proved, 
as is shown in Theorem 1. 
We emphasize that in the 
case $V \le 2$ and Hypothesis B, ${\cal S}$ is not assumed to be 
embedded.  In this context, the requirement that ${\cal S}$ be disk-type 
is necessary; in fact, Wente \cite{WenTub} has given an example of an 
immersed constant mean curvature bridge joining the two faces of a wedge 
and meeting those faces in the contact angles $\gamma_1=\gamma_2=\pi/2$.  
Such a surface cannot lie on a metric sphere.

	In the case $V = 2$, the underlying idea for our work consists 
simply of adjoining a theorem on intersecting surfaces proved by 
Joachimsthal \cite{JoaDem} in 1846, to the method that H. Hopf \cite{HopUbe} 
used in 
1951 to prove that every immersed closed surface of genus zero and 
constant mean curvature is a sphere. Some technical effort will be 
needed to prove our results without superfluous smoothness requirements. 
We present those details in later sections of the paper; in order not 
to obscure the ideas, we outline the proof for the case 
$V\le 2$ in the following section, under some assumptions on smoothness 
of ${\cal S}$ and of certain conformal mappings of ${\cal S}$.  In the 
next following section we provide counterexamples for the case $V= 4$.  
The underlying assumptions used in Section~\ref{(2)} will be justified in 
the ensuing Section 4. 

	In the final Section 5, we present our results for the case 
$V = 3$.  This material is based entirely on comparison procedures, 
and requires an additional hypothesis, although in another sense our 
requirements are weaker (data on a part of the boundary set of 
${\cal Q}$ are allowed 
whereas for $V = 2$ they must be interior to ${\cal Q}$).  It is perhaps 
worth observing that the difficulty in extending the $V = 2$ proof to 
$V = 3$ lies in the singular behavior at the third vertex, of the 
conformal mapping that takes ${\cal S}$ to an infinite strip, with two 
of the vertices going to infinity. That the difficulty is essential 
can be seen from the fact that it is exactly the reason the theorem 
fails when $V\ge 4$. This kind of behavior underscores the importance 
we attach to the material of Section~4, in which the assumptions 
made in Section~2, on asymptotic structure of the mapping at the 
vertices, are justified on the basis of more primitive 
(and seemingly reasonable) hypotheses.

	There is also an essential difficulty extending the $V = 3$ 
proof to $V = 2$, as the additional hypothesis just referred to is 
violated by spheres in the $V = 2$ case, see the comments at the end 
of Section~5.

	We note finally that our method provides as corollary a 
conceptual simplification of a proof given earlier by Nitsche \cite{NitSta}, 
that a disk-type free surface interface of a connected drop, that 
rests in the absence of gravity on a spherical support surface which 
it meets in a constant angle, is necessarily a spherical cap. We 
indicate our improvement to that result in the Appendix to this paper.

\section{Outline of proof: case $V\le 2$.\label{(2)}}
  
We outline here the structure of our proof, assuming to fix the ideas 
that all functions appearing are as  smooth as required by the 
context. It is intuitively clear that the case $V =  1$ does not 
occur; the case $V  = 0$ will be encompassed in the procedure for 
$V  = 2$. In this case the drop can meet only two planes; we may ignore 
any other planes and consider a drop in a dihedral angle, of opening 
$2\alpha$.  We adopt the intersection line ${\cal L}$ as $z$-axis, 
and assume representations $z=u^\pm(x,y)$ for ${\cal S}$ at the two 
vertices ${\cal V}^+$ and ${\cal V}^-$, in a neighborhood of the origin 
interior to the wedge. The condition for ${\cal S}$ to have a tangent 
plane at the origin is precisely the condition that the pair 
$(\gamma_j, \gamma_k)$ lie in the closed rectangle ${\cal Q}$ introduced 
in \cite{ConCap}, see Figure~\ref{rectangle}.  We assume this condition 
satisfied, and assume further that $(\gamma_j, \gamma_k)$ lies interior 
to ${\cal Q}$, so that the tangent plane is not vertical.  The contact 
lines then meet each other  on ${\cal L}$ in the same positive angle 
$2\beta < \pi$ at both vertices, determined entirely by 
$(\gamma_j, \gamma_k)$ and by $\alpha$. 
\begin{figure}[ht]
\centerline{{
         \epsfysize=2.3in
         \epsfxsize=2.3in
         \leavevmode\epsfbox{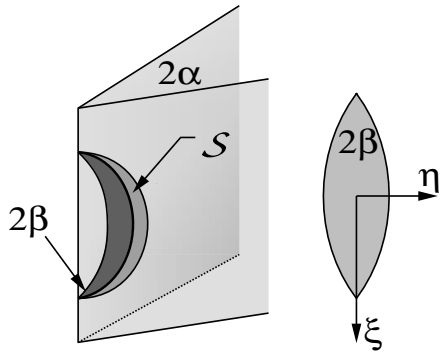}}}
\centerline{ }
\caption{Mapping to lens\label{lens}}
\end{figure}

We map S conformally onto a convex lens (Figure~\ref{lens}) with vertex 
angles $2\beta$ at the points  $P_{-1}= (-1,0)$ and  $P_{+1} = (+1,0)$ in 
the $\zeta = \xi + i\eta$  plane, in such a way that $P_{-1},P_{+1}$ are 
the images of the vertices on ${\cal S}$.  We assume this mapping to 
be sufficiently smooth that the second derivatives of the representation 
$\vec{\bold x}(\xi,\eta)$ of ${\cal S}$  have at worst a singularity 
admitting the estimate 
\begin{equation}\label{estimate}\label{1}	
	|D^2 \vec{\bold x}| = \circ(\rho^{-2})
\end{equation}		
at $P_{-1}$ and $P_{+1}$, $\rho$ being distance to the respective point.  
The mapping
\begin{equation}\label{2}	
	Z = \ln\left({\zeta-1\over{\zeta+1}}\right)
\end{equation}	
takes the lens region to a horizontal infinite strip.  The boundary of 
this domain consists of coordinate lines.  According to a theorem of 
Joachimsthal \cite{JoaDem}, {\em if two surfaces intersect in a constant 
angle, and if the intersection curve is a curvature line on one of the 
surfaces, then it is a curvature line on both surfaces.}  In the present 
case, every curve is a curvature line on the support planes $\Pi_j$;  
since both planes meet ${\cal S}$ in the respective constant angles 
$\gamma_j$, the contact lines are curvature lines on ${\cal S}$. They 
are also coordinate lines in the $Z$ plane;  thus, if we denote by 
$L$, $M$, $N$ the coefficients of the second fundamental form in the 
coordinates $X$, $Y$ of the $Z$-plane, we will have $M = 0$ on the boundary 
of the strip.  Considering the corresponding coefficients $l$, $m$, $n$  
in the $\zeta$ - plane, we note that they can be expressed as scalar 
products of the (assumed continuous) normal vector and the 
derivatives $D^2 \vec{\bold x}(\xi,\eta)$, and thus we have 
by (\ref{estimate}) 
\begin{equation}\label{est2}\label{3}	
	|l|, |m|, |n| = \circ(\rho^{-2}).
\end{equation}
Further, we have from (\ref{2})
\begin{equation}\label{est3}\label{4}	
	\left| {d\zeta\over{dZ}}\right| = \bigcirc(\rho).
\end{equation}

It follows from the Codazzi equations that the expression
\begin{equation}\label{quaddiff}\label{5}	
	\Phi \equiv ((l-n) - 2im)d\zeta^2
\end{equation}
is a holomorphic quadratic differential on ${\cal S}$ 
(see \cite[Chapter 6]{HopDif}).  We thus have
\begin{equation}\label{quaddiff2}\label{6}	
	((L-N) - 2iM) =  ((l-n) - 2im)\left( {d\zeta\over{dZ}}\right)^2
\end{equation}
from which we conclude from (\ref{est2}) and (\ref{est3}) that  
$M\to 0$ uniformly at each end of the strip.  Since   $M$ is harmonic 
in the strip and vanishes on the entire finite boundary, it 
follows from the maximum principle that $M\equiv 0$, and thus that 
its conjugate $(L-N)$ is identically constant.  But $(L-N)\to 0$ at the 
ends for the same reason that $M$ does, and hence $(L-N)\equiv 0$.  
We conclude that ${\cal S}$ is totally umbilic, and must therefore 
be part of a metric sphere in ${\bbR}^3$, as was to be proved. $\Box$

\section{The case $V= 4$\label{(3)}}  

	It is proved in \cite[Sec. 6.4]{FinEqu} that if 
$\pi/4 < \gamma < \pi/2$, then there exists a solution surface 
${\cal S}: u(x,y)$ of the nonparametric constant mean curvature equation
\begin{equation}\label{fineq}	
	{\rm div}\, Tu = 2 {a+b\over{ab}}\cos\gamma
\end{equation}
in a rectangle ${\cal R}$ of arbitrary side lengths $a$ and $b$, such 
that the solution surface meets  all four vertical walls over the sides 
in the constant contact angle $\gamma$.  The solution is uniquely 
determined up to an additive constant.  There are four vertices on 
${\cal S}$, on the vertical lines through the four vertices of ${\cal R}$.  
If $a = b$, then ${\cal S}$ is known explicitly as a lower spherical 
cap.  However, if $a\ne b$, although the solution continues to exist, 
it cannot be spherical.

	Henry Wente pointed out to us that if $\gamma$ is allowed to 
differ on adjacent walls, then an example can be given explicitly.  
Choose $\gamma = 0$ on the two opposite walls of length $a$, and 
$\gamma = \pi/2$ on the two other walls.  Then the lower half of a 
horizontal cylinder of radius $b/2$  provides an explicit surface of 
constant mean curvature $1/b$, meeting the walls in the respective 
angles indicated and having four vertices. 

\medskip

	We devote most of the remainder of this paper to justifying the 
hypotheses we introduced in Section~\ref{(2)} above. In the final section we 
will discuss the case of three vertices. 

\section{The drop configuration; $V\le 2$\label{(4)}} 

	\subsection{Preliminary lemmas} 
We consider a connected drop supported by a finite number of planes, 
with free surface ${\cal S}$ topologically a disk.  We suppose ${\cal S}$ 
to have constant mean curvature and to be differentiable up to the 
(interiors of the) contact lines, where it cuts the planes $\Pi_j$ 
transversally in the respective angles $\gamma_j$ interior to the drop.  
We will make two kinds of hypotheses as to the behavior of ${\cal S}$ 
near the intersection lines; both of them will lead to identical 
further conditions, under which the drop must be spherical.  The 
first of them is relevant only to the case of a drop in a dihedral 
angle formed by two planes $\Pi_1,\Pi_2$, with two vertices on the 
single intersection line ${\cal L}$ (edge of the wedge). We may suppose 
that ${\cal L}$ is oriented vertically.

\begin{lemma}\label{4.1} 
Under the conditions just stated, assume Hypothesis A with respect to 
two distinct points ${\cal V}^\pm$ on ${\cal L}$.  Then ${\cal S}$ is 
symmetric about a horizontal plane $\Pi$, and each half of ${\cal S}$ 
is a graph over  $\Pi$.
\end{lemma}
\noindent{\bf Proof:}  We use the planar reflection method, as introduced 
by Alexandrov \cite{AleUni} and developed for drops on planar surfaces 
by Wente \cite{WenSym}. Let ${\cal V}_0$  be the point on ${\cal L}$ 
midway between ${\cal V}^+$ and ${\cal V}^-$ below it.  We start with a 
horizontal plane $\Pi^+$ as indicated in Hypothesis A.  This plane clearly 
separates ${\cal V}_0$ and ${\cal V}^+$.  We lower the plane 
continuously, reflecting in it the part of ${\cal S}$ that lies above it.  
If a point of tangency of the reflected with the original surface is 
attained, one can show as in \cite{McCSym} via appropriate versions of 
the maximum principle that the lower surface is a reflection of the 
upper one.  Thus the procedure can be continued until the plane 
reaches ${\cal V}_0$.  If no point of tangency (other than at 
${\cal V}_0$)  has been reached till then, it follows that on the 
intersection of the plane $\Pi_0$ through ${\cal V}_0$ with ${\cal S}$, 
the derivative with respect to height of each horizontal distance 
from ${\cal L}$ to the intersection curve ${\cal C}$ of a generic $\Pi$ 
with ${\cal S}$ must be negative.  But if that situation occurs then the 
same procedure, starting with a plane $\Pi^-$ and moving upward, 
would have to yield an earlier point of tangency.  Thus ${\cal S}$ 
is symmetric under reflection in $\Pi_0$.  If the upper and lower parts 
of ${\cal S}$ were not graphs over $\Pi_0$ then a point of tangency 
would be obtained during the procedure, which is not possible since  
${\cal V}^+$ must reflect onto ${\cal V}^-$. $\Box$

\medskip

\begin{figure}[ht]
\centerline{{ 	\epsfysize=2.7in
         	\epsfxsize=2.7in
         \leavevmode\epsfbox{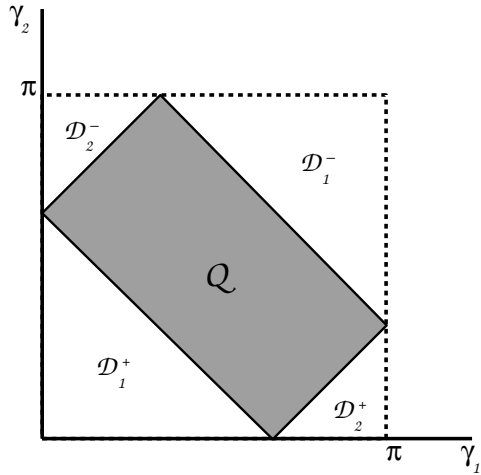}}}
\centerline{ }
\caption{Limit configurations for given data\label{rectangle}}
\end{figure}

 	Lemma~\ref{4.1} reduces the further discussion to the following 
case, which leads independently to conditions under which the existence 
of the vertices ${\cal V}^\pm$ as uniquely determined points can be 
proved. In what follows we continue to assume the smoothness of 
${\cal S}$ in a deleted neighborhood of each intersection line 
${\cal L}$, however ${\cal S}$ is no longer assumed to be embedded, 
or in any sense defined on ${\cal L}$, nor is any growth hypothesis 
introduced with regard to behavior of ${\cal S}$ near ${\cal L}$. 
\begin{lemma}\label{4.2} 
Suppose there is a  neighborhood ${\cal N}$ on a plane  $\Pi$ orthogonal 
to ${\cal L}$ as in Hypothesis B, such that a subset 
${\cal S}^*\subset{\cal S}$ appears as a graph $u^*(x,y)$ over 
${\cal N}\backslash({\cal N}\cap{\cal L})$.  Then $u^*(x,y)$  is bounded in 
${\cal N}\backslash({\cal N}\cap{\cal L})$. 
\end{lemma}
\noindent{\bf Proof:}  The set of all admissible boundary data 
$(\gamma_j,\gamma_k)$, constant on each side of the wedge domain 
determined by the intersecting planes  $\Pi_1,\Pi_2$, can be restricted 
to a square of side length $\pi$.  In \cite{ConCap} it was shown that 
the set of all data that can lead to constant mean curvature graphs 
over such an ${\cal N}$ with tangent planes at ${\cal L}$ lie in an 
inscribed rectangle ${\cal Q}$. The complement of ${\cal Q}$ in the 
square consists of two diagonally opposite domains ${\cal D}_1^\pm$ and 
two diagonally opposite domains ${\cal D}_2^\pm$, 
see Figure~\ref{rectangle}. 
In ${\cal D}_1^\pm$ there can be no such graph over ${\cal N}$, 
regardless of growth 
conditions.  It was shown in \cite{FinLoc} that for data in 
${\cal D}_2^\pm$ graphs meeting $\Pi_1,\Pi_2$ in the prescribed angles 
can under some conditions exist, although they cannot admit tangent 
planes at ${\cal L}$, and in \cite{CheCap} 
it is shown that such graphs must be discontinuous at ${\cal L}$ but are 
nevertheless bounded there.  

In the present case, the data 
$(\gamma_j,\gamma_k)$ cannot lie in a ${\cal D}_1^\pm$ domain as in 
that event there could be no constant mean curvature surface as a graph in 
${\cal N}\backslash({\cal N}\cap{\cal L})$ (\cite[Theorem 3]{ConCap}).  
The boundedness for other data follows from Proposition~1 of 
\cite{LanRad}, 
or alternatively from the material of \cite[Section~7]{CheCap}. 
$\Box$

\begin{lemma}\label{4.3} 
Under the hypotheses of Lemma~\ref{4.2}, suppose additionally that the data 
arise from an interior point of ${\cal Q}$. Then $u^*(x,y)$ can be 
defined at ${\cal L}$ (determining a vertex ${\cal V}^*$ as a point on 
${\cal L}$) so as to have uniformly H\"older continuous derivatives 
up to ${\cal L}$, and the data are achieved exactly near ${\cal L}$ 
by the plane  $\Pi$ tangent to  ${\cal S}$ at the vertex. 
\end{lemma}
The result is obtained by adapting procedures used by 
Simon \cite{SimReg} 
and by Tam \cite{TamReg} 
to prove differentiability in the case $\gamma_1 = \gamma_2$, and then 
by adapting the methods of Lieberman \cite{LieHol} 
or of Miersemann \cite{MieBeh} 
to prove the H\"older continuity of the derivatives.  For details of 
the initial step, see \cite{CheCap}.  
The methods extend without essential change to all data interior 
to ${\cal Q}$.  For data on the boundary of ${\cal Q}$ we still have
\begin{lemma}\label{4.3a}  
Under the conditions of Lemma~\ref{4.2}, if the data lie on the 
interior of the segments $\partial{\cal Q}\cap\partial{\cal D}_1^\pm$, 
then $u^*(x,y)$ can be defined at ${\cal L}$ (determining a vertex 
${\cal V}^*$ as a point on ${\cal L}$) so as to be continuous and 
have continuous unit normal vector. 
\end{lemma}
In this case the first derivatives necessarily become infinite as 
${\cal L}$ is approached. 
The proof 
follows from the procedure of Tam \cite{TamReg},  
see \cite{CheCap}.  

\begin{lemma}\label{4.4}  
Under the hypotheses of Lemma~\ref{4.3},  the two contact lines 
${\cal C}_1^*$ and ${\cal C}_2^*$ intersect on ${\cal L}$ in an angle  
$2\beta$, with  $0 < 2\beta < \pi$.  In the (single angle) case 
$\gamma_1=\gamma_2$  there holds additionally  $0 < 2\beta \le 2\alpha$.
\end{lemma} 
\noindent{\bf Proof:}  According to Lemma~\ref{4.3}, ${\cal S}$ has 
a non-vertical tangent plane $T$ at the vertex, meeting the walls in 
the constant angles $(\gamma_1,\gamma_2)$.  The angle made by 
${\cal C}_1^*,{\cal C}_2^*$  with each other is the same as the angle 
between the 
intersection lines of the tangent plane with the wedge planes.  It 
thus suffices to determine the angle between these lines.  
Set $B_1=\cos\gamma_1$, $B_2=\cos\gamma_2$.  A calculation yields
\begin{equation}\label{8}
	\sin^22\beta = {\sin^22\alpha - (B_1^2+B_2^2 + 2B_1B_2\cos{2\alpha})
			\over{(1-B_1^2)(1-B_2^2)}}.
\end{equation}
It was shown in \cite{ConCap} 
that data arise from an interior point of ${\cal Q}$ if and only if 
the numerator in (\ref{8}) is positive, and thus the first assertion 
follows. To prove the second statement, we note that (\ref{8}) can 
be written in the form
\begin{equation}\label{9} 
	\cos{2\beta} = {B_1B_2 + \cos{2\alpha} 
			\over{\sqrt{1-B_1^2}\sqrt{1-B_2^2}}}
\end{equation}
and thus if $B_1 = B_2 = B$ then
\begin{equation}\label{10} 
	\cos{2\beta} = {B^2 + \cos{2\alpha} 
			\over{1-B^2}}.
\end{equation}
The right side of (\ref{10}) is increasing in $B^2$ and reduces to 
$\cos{2\alpha}$  when $B^2 = 0$; the assertion follows. $\Box$

\begin{lemma}\label{4.5}  
Under the hypotheses of Lemma~\ref{4.3}, the second derivatives of 
$u^*$  are H\"older continuous to the wedge walls, and satisfy  
an estimate $|D^2u^*| < Cr^{-\alpha}$ in terms of distance $r$ 
to the vertex, with  $0 < \alpha < 1$.
\end{lemma}
\noindent{\bf Proof:}   We observe first that by adjoining work of 
Siegel \cite{SieHei}  
to that of Ural'tseva \cite{UraSol}  
and of Gerhardt \cite{GerGlo, GerBou},  
we obtain that locally $u^*(x,y)\in C^{2+\epsilon}$ to the wedge walls 
$(\Pi_1\cap\Pi_2)\backslash{\cal V}^*$, with $0 < \epsilon <  1$.  
By Lemma~\ref{4.3}, $u^*(x,y)\in C^{1+\epsilon}$ in the closed wedge 
domain.  We next observe that ${\cal S}$  can be represented locally 
near ${\cal V}^*$ as a graph $u(x,y)$ over its tangent plane.  Since 
$u^*(x,y)\in C^{1+\epsilon}$, there follows $|u| < Cr^{1+\epsilon}$ 
in polar coordinates centered at ${\cal V}^*$.  We have also that 
$u(x,y)$ satisfies a uniformly elliptic equation of the form
\begin{equation}\label{11} 
	a_{ij}(x)u_{x_ix_j} = 2H
\end{equation}		
near ${\cal V}^*$, with coefficients and Dirichlet boundary data 
that are H\"older continuous in the closed corner domain.  
Such solutions were studied by Azzam \cite{AzzBeh}, 
who obtained the stated growth estimate on second derivatives. $\Box$
\begin{lemma}\label{4.6}  
Under the hypotheses of Lemma~\ref{4.3}, a neighborhood of ${\cal V}^*$ 
on ${\cal S}$ can be mapped 1-1conformally onto a corresponding 
neighborhood of a rectilinear angle, of opening $2\beta$, such that 
the inverse representation of ${\cal S}$ over the angular neighborhood 
is locally of class $C^{2+\alpha}$ up to the rectilinear sides, and of 
class $C^{1+\alpha}$ to the vertex image ${\cal V}'$.  For the second 
derivatives of the position vector $\vec{\bold x}$ in this representation, 
there holds $|D^2\vec{\bold x}|<Cr^{\epsilon-1}$, with  $0 < \epsilon < 1$.
\end{lemma} 
\noindent{\bf Proof:}   We consider ${\cal S}$ in local representation 
$u(x,y)$ over its tangent plane $\Pi$ at ${\cal V}^*$.  In view of 
Lemma~\ref{4.3}, the quantities
\begin{equation}\label{12} 
	E= 1 + u_x^2,\ F=u_xu_y,\  G = 1+u_y^2
\end{equation}	
are defined in a wedge region ${\cal W}$ determined by the projections of 
${\cal C}_1^*$, ${\cal C}_2^*$ onto $\Pi$;  in view of Lemma~\ref{4.5}, 
they are in 
class $C^\alpha$ on ${\cal W}$, with $E = G = 1$, $F = 0$ at ${\cal V}^*$, 
and in class $C^{1+\alpha}$ in ${\cal W}\backslash{\cal V}^*$.     
\begin{figure}[ht]   
\centerline{{ 	\epsfysize=2.7in
         	\epsfxsize=2.7in
         \leavevmode\epsfbox{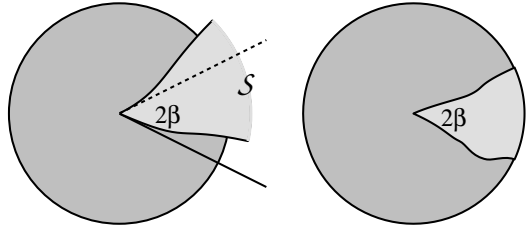}}}
\centerline{ }
\caption{Local mapping at vertex\label{local}}
\end{figure}
We extend these functions to functions with the same smoothness 
properties defined in a (small) disk about ${\cal V}^*$ 
(Figure~\ref{local}), and observe that non-singular solutions 
$\zeta \in {\Bbb C}$ of the {\em Beltrami Equations} 
\begin{equation}\label{13}
	{d\zeta\over{d\bar{z}}} = \lambda {d\zeta\over{dz}}
\end{equation}
where
$$
	\lambda = {E-G+2iF\over{E+G+2\sqrt{EG-F^2}}}
$$
in the disk determine conformal maps of the portion of ${\cal S}$ 
that projects onto ${\cal W}$.  Since $E$, $F$, $G$  are H\"older continuous, 
there exists a local solution $\zeta = \xi + i\eta$ about ${\cal V}^*$ 
with H\"older continuous derivatives and non-vanishing Jacobian 
determinant \cite{CouMet},  
mapping ${\cal W}\leftrightarrow {\cal W}'$.  Since at ${\cal V}^*$, 
$E = G = 1$, $F = 0$, the mapping is conformal between the planar 
domains at ${\cal V}^*$ and 
thus the vertex angle $2\beta$ remains unchanged in the ${\cal W}'$ 
coordinates. We may assume the angles oriented as in Figure~\ref{local}.  

\begin{figure}[ht]   
\centerline{{ 	\epsfysize=3.0in
         	\epsfxsize=3.0in
         \leavevmode\epsfbox{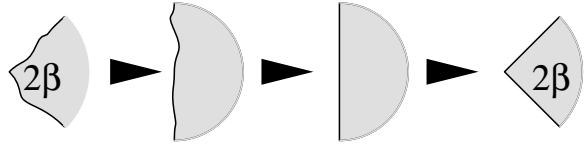}}}
\centerline{ }
\caption{Successive mappings\label{maps}}
\end{figure}

The mapping $\Xi = \zeta^{\pi/2\beta}$ opens the wedge to a domain 
bounded locally by a H\"older differentiable curve ${\cal C}$.  A further 
mapping $Z = F(\Xi)$, again H\"older differentiable and invertible 
to the image of ${\cal V}^*$, takes ${\cal C}$ onto a linear segment 
(Figure~\ref{maps}).  
Finally, $\Psi = Z^{2\beta/\pi}$ completes the mapping to a rectilinear 
wedge domain of opening $2\beta$.

	According to Lemma~\ref{4.3}, the surface representation 
$\vec{\bold x}(x,y)$ is H\"older differentiable to ${\cal V}^*$.  Since 
the mapping $\zeta(x,y)$ is H\"older differentiable to ${\cal V}^*$, 
it follows that $\vec{\bold x}$ has the same property in the $(\xi,\eta)$ 
variables.  The asserted properties of the first derivatives of 
$\vec{\bold x}$ in the $\Psi$ variables follow by tracing through 
the mappings, each of which is either smooth at the vertex or is 
an explicitly known power mapping.  The singularities in the power 
mappings cancel each other.

To obtain the growth estimate on second derivatives, we observe that 
in the $\Psi$-variables, the representation $\vec{\bold x}$ satisfies 
the equation
\begin{equation}\label{14} 
	\Delta \vec{\bold x} = 2 H \vec{N}
\end{equation}		 
where $\vec{N}$ is a unit normal to ${\cal S}$.  By the material 
just proved, $\vec{N}$ is H\"older continuous to the vertex, while 
$\vec{\bold x}$ is in $C^{2+\alpha}$ locally to the wedge walls, and in 
$C^{1+\epsilon}$ to the vertex. The result then follows as in the 
second part of Lemma~\ref{4.5}. $\Box$
\begin{lemma}\label{4.7}  
Under the hypotheses of Lemma~\ref{4.3}, every  conformal map of a 
neighborhood of ${\cal V}^*$ on ${\cal S}$ onto a rectilinear wedge 
neighborhood of opening $2\beta$, with vertex going to vertex, 
leads to a representation for ${\cal S}$ with the smoothness 
properties described in Lemma~\ref{4.6}.
\end{lemma}
\noindent{\bf Proof:}  A given mapping onto the $\hat{\zeta}$ plane, in 
conjunction with the particular mapping 
$\hat{\Xi} = \hat{\zeta}^{\pi/2\beta}$ applied in the proof of 
Lemma~\ref{4.6}, leads to a conformal mapping into itself of a half disk 
of which the image of a diameter lies on a diameter, with origin $O$  
going into itself.  This map can be extended by reflection into a 
univalent conformal map of a disk containing $O$, and is therefore 
analytic with non-zero Jacobian at $O$. The mappings to the wedge 
are effected by an identical mapping for the two functions, with the 
requisite smoothness properties, and hence the H\"older 
differentiability of the constructed function leads to the same 
property for the given one.  Similarly, the second derivatives of the 
two functions have the same H\"older growth exponent at the respective 
images of $O$. $\Box$

\smallskip

	As a consequence of Lemma~\ref{4.7}, we obtain immediately, 
since $l$, $m$, $n$ are scalar products of the unit normal to 
${\cal S}$ with second derivatives of the position vector:
\begin{corollary}\label{c4.7}  
In terms of the conformal parameters introduced in Lemma~\ref{4.6}, 
the coefficients $l$, $m$, $n$ of the second fundamental form satisfy
\begin{equation}\label{15}
	\sqrt{l^2 + m^2 + n^2} < C r^{\epsilon -1}
\end{equation}
near the vertex image. 
\end{corollary}


\subsection{Main theorems, $V \le 2$}
\begin{theorem}\label{t1}  
Under either of the hypotheses A or B, if  $V = 2$ and ${\cal S}$ is 
topologically  a  disk, 
and if the contact angles 
$\gamma_1$, $\gamma_2$ lie interior to the rectangle ${\cal Q}$, 
then ${\cal S}$ is metrically spherical. 
\end{theorem}
\noindent{\bf Proof:}  If Hypothesis A holds, we conclude by 
Lemmas \ref{4.1} to \ref{4.4} that ${\cal S}$ is symmetric about a 
plane $\Pi$ orthogonal to the intersection line ${\cal L}$, can be 
represented globally by functions with H\"older continuous first 
derivatives, and forms at each vertex an angle $2\beta$ given by 
(\ref{10}).  We can therefore map ${\cal S}$ conformally onto a lens 
domain (Figure~\ref{lens}) bounded by circular arcs meeting at angle 
$2\beta$, with vertices going into vertices at the points $\zeta = \pm 1$.  
This configuration is locally related to a rectilinear angle via a 
linear fractional transformation; we may thus conclude from 
Corollary~\ref{c4.7} that in the lens coordinate $\zeta$, 
the coefficients $l$, $m$, $n$  satisfy (\ref{15}).  We now apply 
the mapping (\ref{2}) and the invariance of the form $\Phi$ as in 
Section~\ref{(2)}, arriving at the desired conclusion by 
the identical reasoning. 

	If Hypothesis B holds, we are assured directly of the 
hypotheses of Lemma~{4.2}; the remainder of the reasoning then 
proceeds without change. $\Box$

\medskip

	If we assume additionally the Vertex Condition, then we 
can provide a more inclusive formulation of Theorem~\ref{t1}.  We 
observe that the variational condition characterizing the mechanical 
equilibrium of the drop surfaces is not affected by the presence or 
absence of support surfaces that do not meet ${\cal S}$.  Thus, 
{\em if a property of ${\cal S}$ has been determined by its interaction 
with certain support planes which it contacts, the removal of planes 
which do not contact ${\cal S}$ will not affect that result.}\/  Using 
this observation, we are able to characterize equivalence classes 
of configurations in terms of a few particular cases.  Specifically, 
we find:
\begin{lemma}\label{4.8}  
Assume the Vertex Condition and that one of the Hypotheses A or B holds 
at each vertex, with data arising from the interior of ${\cal Q}$ 
or from  the interior of $\partial{\cal Q}\cap \partial{\cal D}_1^\pm$.  
Then if $V = 0$, every drop configuration is either a closed surface 
without boundary or else it can be realized by a drop on a single 
plane.  The case $V = 1$ does not occur.  If $V = 2$ then the 
configuration can be realized by a drop in a dihedral angle 
(wedge).  If $V = 3$ then the configuration is equivalent either 
to a drop covering the vertex in a trihedral angle, or else to a 
drop covering the (planar) base of a cylindrical container, 
whose side walls consist of three planes, no two of which are 
parallel. 
\end{lemma}
\noindent{\bf Proof:}  We have assumed the boundary ${\cal C}$ of 
${\cal S}$ to be locally smooth on each face;  by 
Lemmas~\ref{4.3}, \ref{4.3a} and \ref{4.4} it is piecewise smooth at 
each vertex, that is, continuous with a jump in unit tangent vector.  
Since the entire configuration is compact and since ${\cal S}$ is 
disk type, ${\cal C}$ is globally a 
piecewise smooth closed curve (which may conceivably have 
self-intersections).  We may 
start with any point on ${\cal C}$, and 
traverse 
${\cal C}$ in any chosen direction. 

If $V = 0$ and ${\cal C}$ is the null set, then ${\cal S}$ is closed 
and without boundary.  If $V=0$ and ${\cal C}$ is non-null, we start with an 
arbitrary point of ${\cal C}$, which will be on a support plane $\Pi$, 
and  
traverse 
${\cal C}$ in one of the two possible directions.  ${\cal C}$ cannot 
enter a plane distinct from $\Pi$ across an intersection line with $\Pi$, 
as, by the Vertex Condition, that would create a vertex.  
Thus all planes distinct from 
$\Pi$ can be removed without affecting the variational conditions 
determining ${\cal S}$.  We are left with a drop with disk-type free 
surface resting on a single plane.
 
	Suppose $V = 1$.  We 
traverse 
${\cal C}$ 
beginning at its single vertex ${\cal V}$, along one of 
the two intersecting planes $\Pi_1$ and $\Pi_2$.  ${\cal C}$ is a closed 
curve and on each of $\Pi_1$ and $\Pi_2$  it 
contains points distinct from the intersection line ${\cal L}$.  
${\cal C}$ cannot meet other planes, as by the 
Vertex Condition that would create a new vertex. Therefore it cannot close 
without crossing again over ${\cal L}$ at a point distinct from 
${\cal V}$, which cannot occur by the Vertex Condition.  We conclude that 
this case does not occur. 

	The same reasoning shows that if $V = 2$, the second vertex 
must be a (distinct) point on the same intersection line ${\cal L}$.  Again no 
other planes can be contacted, and we may thus discard all planes 
distinct from $\Pi_1$ and $\Pi_2$.

	Suppose finally that $V = 3$.  Starting with a given vertex 
${\cal V}_{12}$ on the intersection line ${\cal L}_{12}$ joining 
planes $\Pi_1$ and $\Pi_2$, we follow 
${\cal C}$ from $\Pi_1$ across ${\cal L}_{12}$ onto $\Pi_2$ until the 
vertex  ${\cal V}_{23}$ 
appears on the intersection line  ${\cal L}_{23}$  between $\Pi_2$ and 
$\Pi_3$.  We claim that ${\cal L}_{23}$ is distinct from ${\cal L}_{12}$.  
Otherwise, since no three planes can intersect along ${\cal L}_{12}$, 
${\cal C}$ must continue back onto $\Pi_1$ and then to the third vertex 
${\cal V}_{31}$ on an intersection line ${\cal L}_{31}$.  If ${\cal L}_{31}$ 
coincides again with ${\cal L}_{12}$, then ${\cal C}$ must continue back 
onto $\Pi_2$ and could not join the initially chosen points on $\Pi_1$ 
without crossing a fourth vertex.  Also, if ${\cal L}_{31}$ is distinct from 
${\cal L}_{12}$ the same contradiction arises.  

Thus, we may assume that ${\cal L}_{23}$ is distinct from ${\cal L}_{12}$, 
and that ${\cal C}$ crosses ${\cal L}_{23}$ at the vertex ${\cal V}_{23}$ 
onto a plane $\Pi_3$ distinct from $\Pi_1$ and $\Pi_2$ (and proceeds to the 
third vertex ${\cal V}_{31}$ on an intersection line ${\cal L}_{31}$ 
distinct from ${\cal L}_{12}$ and ${\cal L}_{23}$, by the same argument).  

We observe that ${\cal V}_{12}$   and ${\cal V}_{23}$  share $\Pi_2$  
as a common plane serving as one of the intersecting planes for both 
vertices.  Similarly $\Pi_3$ is shared by ${\cal V}_{23}$ and 
${\cal V}_{31}$. 

	We assert that  $\Pi_1$  is shared by  ${\cal V}_{31}$ and by  
${\cal V}_{12}$.  For by construction, it is one of the sides for 
${\cal V}_{12}$.  If ${\cal C}$  were to continue 
through ${\cal V}_{31}$ onto a plane distinct from $\Pi_1$, it 
would have to pass through still another vertex before returning to 
${\cal V}_{12}$, contrary to hypothesis.  Thus, ${\cal C}$ encounters 
only the three planes $\Pi_1$, $\Pi_2$, $\Pi_3$ of the supporting 
family, of which no two can be parallel, and we conclude also that 
${\cal S}$ can encounter only those planes, as otherwise it would not 
be disk-type.  All other planes can be deleted without affecting the 
variational conditions or the configuration.

	If $\Pi_1$, $\Pi_2$, $\Pi_3$ share a common point ${\cal O}$, 
then we have a drop covering the vertex of a trihedral angle 
(see Figure 6).  The other possibility is that the normals of         
$\Pi_1$, $\Pi_2$, $\Pi_3$ lie in a common plane, or equivalently 
that the three planes share a common (generating) direction.  In this case 
we can replace all other planes of the supporting family with a single 
plane, situated far enough along the generating direction so as not to 
meet ${\cal S}$.  We are done. $\Box$

\medskip

We are now prepared to prove:
\begin{theorem}\label{t2}  
Suppose ${\cal S}$ has $V \le 2$ vertices, that the Vertex Condition 
holds, and that the data on interior 
points of any two adjacent support planes come from interior points of 
${\cal Q}$.  We assume further either that  the hypotheses of 
Lemma~\ref{4.1} are fulfilled, or else that the hypotheses of 
Lemma~\ref{4.2} hold, with respect to each vertex.  Then  ${\cal S}$  
is metrically spherical. 
\end{theorem}
\noindent{\bf Proof:}  We use Lemma~\ref{4.8}.  
If $V =0$ then the configuration is equivalent to a drop on a single plane, 
which it meets in a constant contact angle $\gamma$,  $0 < \gamma < \pi$.  
If the contact set ${\cal C}$ is null, then the statement is equivalent 
to Hopf's theorem \cite{HopUbe}.  If ${\cal C}$ is non-null, then in 
suitable local coordinates near any of its points, 
${\cal S}$ can be represented as a graph $u(x,y)$ meeting a vertical 
wall in angle $\gamma$.  As in the proof of Lemma~\ref{4.5} above, 
we find that $u$ is twice H\"older differentiable to the boundary.  
It follows that ${\cal S}$ can be mapped conformally to the unit disk 
$|\zeta|<1$, with $l$, $m$, $n$  H\"older continuous to the boundary.  
As noted above, $\Phi \equiv ((l-n) - 2im)\,d\zeta^2$ is a 
holomorphic quadratic differential in the disk. 

	The mapping (\ref{2}) takes the disk (and hence ${\cal S}$) 
to an infinite horizontal strip, so that the contact line ${\cal C}$ goes 
into the two bounding coordinate lines, with $l, m, n \to L,M,N$.  By 
Joachimsthal's theorem \cite{JoaDem}  
these lines are curvature lines on ${\cal S}$, so that 
$M = 0$ on them.  Following the discussion in Section~\ref{(2)}, 
we see that $L$, $M$, $N$ all tend uniformly to zero at the ends of the 
strip.  Since $M$ and $L - N$  are harmonic in the strip, we 
conclude from the maximum principle first that $M\equiv 0$, and then 
that $L\equiv N$.  ${\cal S}$ is therefore totally umbilic and hence 
must be spherical. 


\medskip 

The case $V=1$ is vacuous by Lemma~\ref{4.8}, and $V=2$ reduces to a 
dihedral angle, which case is covered by Theorem~\ref{t1}. $\Box$

\section{The Case $V = 3$}\label{(5)}
\begin{figure}[ht]   
\centerline{{ 	\epsfysize=3.0in
         	\epsfxsize=3.0in
         \leavevmode\epsfbox{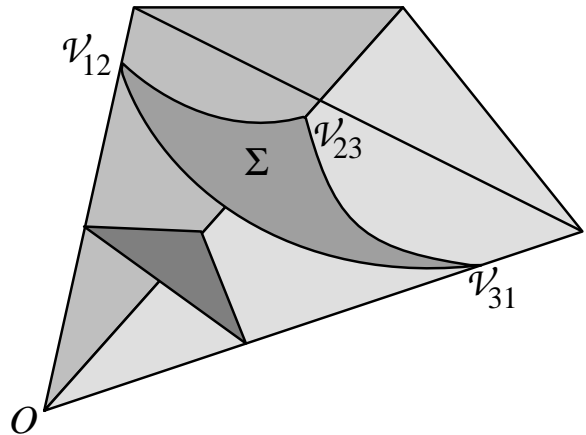}}}
\centerline{ }
\caption{Equivalent spherical cap configuration; the dark triangle 
indicates a conceivable support plane, which would not affect the shape 
of $\Sigma$.  The spherical cap solution $\Sigma$ shown is by 
Theorem~\ref{t3} the unique solution.\label{equiv}}
\end{figure}

By Lemma~\ref{4.8}, we may assume the configuration to be a trihedral 
angle bounded by three planes $\Pi_1$, $\Pi_2$, $\Pi_3$, with the angle 
vertex $O$ covered by the liquid, or else a cylinder with sides 
$\Pi_1$, $\Pi_2$, $\Pi_3$ and closed at one end by a base (as illustrated 
conceptually in Figure~\ref{equiv}).

We consider first the trihedral angle.  In this configuration an explicit 
spherical cap solution $\Sigma$ 
can be given (Figure~\ref{equiv}) corresponding to any prescribed 
mean curvature $|H|>0$, and contact angles $\gamma_1$, $\gamma_2$, $\gamma_3$, 
each pair of which lies in ${\cal Q}$.  Specifically,  the condition 
that $\Sigma$ cuts $\Pi_1$ in angle $\gamma_1$ and $\Pi_2$ in angle 
$\gamma_2$ and encloses a segment of the edge ${\cal L}_{12}$ 
is precisely that $(\gamma_1,\gamma_2)\in{\cal Q}$.  The centers of all 
such spheres of radius $1/|H|$ lie on a line ${\cal L}_{12}'$ parallel to 
${\cal L}_{12}$.  (For a proof that ${\cal L}_{12}'$ is uniquely 
determined, see \cite{FinLoc}).  
Similarly, the condition that $\Sigma$ cuts $\Pi_2$ in angle $\gamma_2$ 
and $\Pi_3$ in angle $\gamma_3$ and encloses a portion of 
${\cal L}_{23}$ is that $(\gamma_2,\gamma_3)\in{\cal Q}$.  The line  
${\cal L}_{23}'$ of centers of such spheres is not parallel to 
${\cal L}_{12}'$, and 
it lies in a plane of centers of spheres of the given radius that meet 
$\Pi_2$ in angle $\gamma_2$.  This plane contains ${\cal L}_{12}'$.  
Thus ${\cal L}_{23}'$ intersects ${\cal L}_{12}'$ in a unique point 
$P_{123}$, which provides the center of a sphere of radius $1/|H|$ 
with the stated property.  We note that the procedure provides exactly 
two spherical caps, each part of the same sphere, for which the mean 
curvature vectors are directed 
respectively into or out of the region cut off by the planes and the 
spherical surface $\Sigma$.  A corresponding sphere of any 
radius meeting the walls 
in the same angles can then be obtained by a similarity transformation.  
An examination of the procedure shows that it applies also for data in 
the interior of $\partial{\cal Q}\cap \partial{\cal D}_1^\pm$.  
In that case the sphere will meet the corresponding 
line or lines in a single point, rather than enclose a segment as above.  
On the other hand, if any data pair lie exterior to the sets considered,
then we find from the results of \cite{ConCap}  
that no spherical solution can exist.  We conclude for trihedral angles 
that {\em if $V = 3$, 
then every surface of mean curvature $|H|>0$ that meets the bounding planes 
in the given angles is spherical if and only if any spherical solution 
is uniquely determined among solutions of the same mean curvature. }

\smallskip

We now introduce a condition under which the spherical surface is the 
only possibility.  We model our discussion on the procedure used by 
Vogel in \cite{VogUni},  
who restricted attention to embedded surfaces bounded on a {\em smooth} 
supporting cone, in the sense that the support surface is assumed to 
be conical and to cut a unit sphere centered at the vertex in a curve 
with continuously turning tangent and lying in a hemisphere.  Under the 
convention that $H > 0$ if the mean curvature vector points 
exterior to the region bounded between ${\cal S}$ and the vertex, he 
was able to prove that if $H > 0$ then ${\cal S}$ is uniquely 
determined by the boundary angle, which he assumed constant.  His proof 
as given in \cite{VogUni}  
does not apply to the (non-smooth) configurations and discontinuous 
boundary angles studied in the present paper.  However, we are able to 
extend the result.  
\begin{theorem}\label{t3}  
Assume the hypotheses of Lemma~\ref{4.2}, that $V = 3$, and that the 
Vertex Condition holds.  Given data that arise from the interior of 
${\cal Q}$ or from the interior of $\partial{\cal Q} \cap 
\partial{\cal D}_1^\pm$, if the support surface given by 
Lemma~\ref{4.8} is a trihedral angle, 
then for any constant $H\ge 0$ any embedded disk-type surface ${\cal S}$ of 
mean curvature $H$ that lies interior to the angle and meets the three 
plane pairs, in whose intersections the vertices lie, in the prescribed 
angles, is metrically spherical. 
\end{theorem}
\noindent{\bf Proof:}  As noted just above, if $H>0$ then it suffices 
to show the uniqueness of a spherical cap solution in a trihedral angle as 
in Figure~\ref{equiv}.  
More generally, we will prove the uniqueness 
of any surface of the type considered.  By Lemma~\ref{4.2}, if we choose 
any of the three intersection lines ${\cal L}$ as vertical axis in a 
Euclidean frame, then the height $u(x)$ of any such ${\cal S}$ is 
bounded near ${\cal L}$.  By Tam's theorem \cite{TamReg},  
see also \cite{CheCap}, 
$u(x)$ has first derivatives continuous to ${\cal L}$, and hence is itself 
continuous in the closure of a neighborhood ${\cal N}$ in a base plane.  Thus, 
the closure of ${\cal S}$ can be represented by continuous functions. 

\begin{figure}[ht]   
\centerline{{ 	\epsfysize=3.0in
         	\epsfxsize=3.0in
         \leavevmode\epsfbox{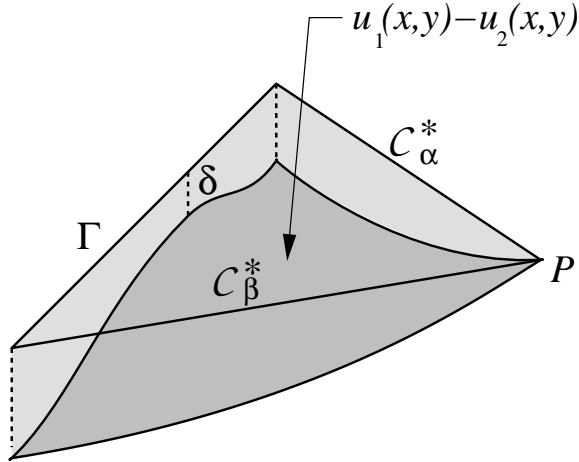}}}
\centerline{ }
\caption{Construction for comparison proof\label{diff}}
\end{figure}

We follow in outline Vogel's reasoning.  
If there were two distinct surfaces ${\cal S}_1$ and ${\cal S}_2$ with 
the same $H$ and the same boundary conditions, there would be a point 
on one of the surfaces, say ${\cal S}_1$, that is exterior to the region 
${\cal I}_2$ bounded between ${\cal S}_2$ and $O$ 
(see Figure~\ref{equiv}). 
Now scale ${\cal S}_1$ by a factor $\lambda < 1$, so that the scaled 
surface $\lambda{\cal S}_1$ lies in the closure of ${\cal I}_2$, and 
that there is at least one point on the closure of $\lambda{\cal S}_1$ 
that lies on the closure of ${\cal S}_2$.  If any such point lies at a 
point of that closure distinct from the vertices, we can proceed as in 
\cite{VogUni}  
and derive a contradiction from a touching principle.  We may thus assume 
that all contact points lie at the vertices.  Let ${\cal V}$ denote 
such a vertex, and ${\cal L}$ the corresponding intersection line, 
which we adopt as $z$-axis.  We may then introduce a segment $\Gamma$ 
in ${\cal N}$ as in Figure~\ref{diff},  
cutting off with the projections ${\cal C}_\alpha^*$, ${\cal C}_\beta^*$  
of the contact lines through ${\cal V}$ a closed triangle $T$ over which 
$u_1(x,y) \le u_2(x,y)$, equality holding only at the  projection $P$ of 
${\cal V}$. 

	On $\Gamma$, $u_1(x,y) \le u_2(x,y)-\delta$, with $\delta > 0$.  
Also, $\lambda{\cal S}_1$ and ${\cal S}_2$ meet the vertical planes 
over ${\cal C}_\alpha^*$, ${\cal C}_\beta^*$ in identical angles 
$\gamma_\alpha$, 
$\gamma_\beta$.  In $T$, the functions $u_2(x,y)$ and 
$u_1^\delta(x,y)\equiv u_1(x,y) + \delta$ satisfy respectively the equations
\begin{eqnarray}\label{(16)}	
	{\rm div}\, Tu & = & 2H \\
	{\rm div}\, Tu & = & {2\over{\lambda}}H > 2H, 
\end{eqnarray}
where $Tu\equiv \nabla{u}/\sqrt{1+|\nabla{u}|^2}$.
	
On $\Gamma$, $u_2\ge u_1^\delta$.  On ${\cal C}_\alpha^*$, 
${\cal C}_\beta^*$ we have  
$\nu\cdot Tu_2 = \ {\rm cosine\ of\ boundary\ angle}\  = 
\nu\cdot Tu_1^\delta$.  
The point $P$ is a set of linear Hausdorff measure zero.  By the 
comparison principle Theorem~5.1 of \cite{FinEqu},  
there follows $u_2\ge u_1^\delta$ throughout $T$.  But by the 
construction, $u_1^\delta > u_2$ near $P$.  This contradiction establishes 
the theorem when $H>0$. 

	If $H = 0$ the discussion requires some changes in detail; 
notably the configuration is no longer uniquely determined, as scaling 
of any given solution leads to a continuum of further solutions.  
Also a planar solution is determined by data on only two of the three 
support planes; the data on the third plane lead to an overdetermined 
problem.  On the other hand, scaling does not affect the curvature, 
thus permitting freedom in the direction in which the surface is 
scaled.  We proceed as follows:

	We suppose given an embedded surface ${\cal S}$ of 
mean curvature $H = 0$ in the trihedral angle formed by the planes 
$\Pi_1$, $\Pi_2$, $\Pi_3$, meeting the planes in the angles 
$\gamma_1$, $\gamma_2$, $\gamma_3$, and having tangent planes continuous 
to the three vertices.  Denote by $\Pi$ a plane tangent to ${\cal S}$ 
at the vertex ${\cal V}_{12}$.  We intend to show that $\Pi$ and 
${\cal S}$ are identical. 

	According to the construction, $\Pi$ meets the planes 
$\Pi_1$ and $\Pi_2$ in the angles $\gamma_1$, $\gamma_2$.  It meets 
$\Pi_3$ in a constant angle $\hat{\gamma}_3$ (if $\Pi$ is parallel 
to $\Pi_3$ then $\hat{\gamma}_3 = 0$), and we suppose initially 
that $\hat{\gamma}_3 \le\gamma_3$.  We move $\Pi$ rigidly away from 
$O$,keeping its unit normal vector constant, until ${\cal S}$ is 
contained strictly in the region bounded 
between $\Pi$ and $O$, and then move $\Pi$ back toward $O$ until a first 
point of contact appears. 

If $\hat{\gamma}_3 < \gamma_3$, then such a point cannot appear on 
$\Pi_3$.  In this case, the reasoning of Vogel excludes the appearance 
of other contact points distinct from the vertices unless 
${\cal S}$ and $\Pi$ coincide.  If $\hat{\gamma}_3 = \gamma_3$, then 
Vogel's reasoning shows directly that any such contact point must be a 
vertex.

	We may thus suppose as before that all initial contact points 
occur at vertices.  If such a point appears at ${\cal V}_{12}$ then 
we may complete the reasoning as above.  We therefore suppose an initial 
contact point at one of the other vertices, which we denote by ${\cal V}$; 
this implies in particular that $\Pi$ is not parallel to $\Pi_3$, 
and thus that we can adopt the intersection line ${\cal L}$ 
through ${\cal V}$ as $z$-axis, with local representations 
$u_1(x,y)$, $u_2(x,y)$ for ${\cal S}$ and for $\Pi$.  The configuration 
is again illustrated by Figure~\ref{diff}, 
however both functions $u_1(x,y)$ and 
$u_2(x,y)$ satisfy the same equation ${\rm div}\, Tu = 0$, and the 
relation $\nu\cdot Tu_2 = \nu\cdot Tu_1^\delta$ is now replaced by 
$\nu\cdot Tu_2 \ge \nu\cdot Tu_1^\delta$.  Since $u_2\ge u_1^\delta$ on 
$\Gamma$, the comparison principle again yields $u_2\ge u_1^\delta$ near 
$P$, a contradiction. 

	There remains the possibility $\hat{\gamma}_3 > \gamma_3$.  
In this event $\Pi$ cannot be parallel to $\Pi_3$, and we can move 
$\Pi$ inwards toward $O$ and then outward till an initial contact point 
with ${\cal S}$ appears.  An analogous reasoning then applies, 
and we conclude finally that ${\cal S}$ and $\Pi$ coincide, so 
that ${\cal S}$ is planar.  $\Box$

\medskip 

	As a particular consequence of the above reasoning, we obtain:
\begin{corollary} \label{cor3}  
If a disk-type minimal surface ${\cal M}$ lies in a trihedral angle, 
meets the sides of that angle in constant angles $\gamma_1$, $\gamma_2$, 
$\gamma_3$ and has a tangent plane continuous to the edges, 
then ${\cal M}$ is a plane. 
\end{corollary}

\smallskip

\noindent{\bf Remark:}  In spherical coordinates, a surface 
$u(\theta,\phi)$ of mean curvature $H$ satisfies the equation
\begin{equation}\label{(17)}
{\partial\over{\partial\theta}}{u_\theta\over{W}} + 
{\partial\over{\partial\phi}}{u_\phi\sin^2\phi\over{W}} = 
	2\left({\sin\phi\over{W}}+H\right)u\sin\phi,
\end{equation}
where $W=\sqrt{(u^2 + u_\phi^2)\sin^2\phi + u_\theta^2}$.
The method of Ambrazevi\u{c}ius \cite{AmbSol,AmbFin}  
applied to such an equation when $H \ge 0$ yields a comparison principle 
analogous to that of \cite[Theorem~5.1]{FinEqu},   
and shows directly the uniqueness of any solution with mixed Dirichlet and 
boundary angle data.  For surfaces that admit such a representation, 
this result could have replaced the procedure we used to prove 
Theorem~\ref{t3}.  The results in \cite{AmbSol,AmbFin}  
have independent interest, and assure the existence of solutions in 
conical regions, with prescribed boundary data and contact angle. $\Box$

\medskip

 
\medskip

In the case of a cylinder, the sign of $H$ is irrelevant:
\begin{theorem}\label{t4}  
Assume the hypotheses of Lemma~\ref{4.2}, that $V = 3$, and that the 
Vertex Condition holds.  Given data that arise from the interior of 
${\cal Q}$ or from the interior of $\partial{\cal Q} \cap 
\partial{\cal D}_1^\pm$, if the support surface given by 
Lemma~\ref{4.8} is a cylinder, 
then any embedded disk-type surface ${\cal S}$ of constant 
mean curvature that meets the three plane pairs, in whose 
intersections the vertices lie, in the prescribed angles, 
is metrically spherical. 
\end{theorem}
\noindent{\bf Proof:}  A discussion analogous to that proceeding 
Theorem~\ref{t3} shows that {\em any angle data of the type described in 
the theorem determines a unique (up to translation) spherical solution 
$\Sigma$}.  Denote the mean curvature of this spherical solution with 
respect to the normal pointing 
out of the enclosed volume by $H_0$. 
 
\begin{figure}[ht]   
\centerline{{ 	\epsfysize=1.5in
         	\epsfxsize=1.5in
         \leavevmode\epsfbox{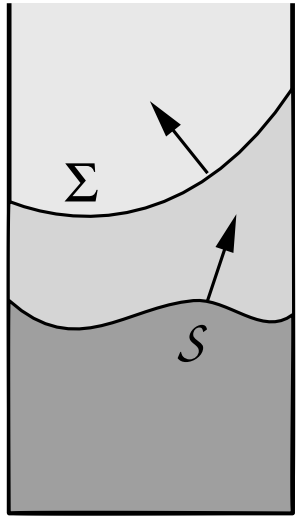}}}
\centerline{ }
\caption{Construction for comparison proof, cylindrical case\label{cyl}}
\end{figure}

Let ${\cal S}$ be any non-spherical solution satisfying the same boundary 
conditions and with constant mean curvature $H$.  If $H\le H_0$, then 
translate ${\cal S}$ along the generating direction (see Figure~\ref{cyl}) 
until the drop volume determined by $\Sigma$ is contained in that determined 
by ${\cal S}$ and yet ${\cal S}\cap \Sigma \ne \phi$.  Proceeding as in 
the proof of Theorem~\ref{t3}, we obtain a contradiction by the comparison 
principle.  If $H\ge H_0$, then we translate $\Sigma$ to obtain the same 
contradiction.  $\Box$

\medskip

While we have extended Theorem~\ref{t3} to the (limiting) case in which 
the three planes have coplanar normals, 
so that the supporting configuration is cylindrical,  we 
note that Theorem~\ref{t4} is false for polyhedral angles of more than three 
sides; in fact, the spherical cap solution is determined by any three 
sides of the angle.  Thus, given any three such sides, the remaining 
configuration is uniquely determined by the requirement that a 
spherical cap be a solution. 

	It should be noted again that our demonstrations for the 
case $V = 3$ differ in essential ways from those we present when $V = 2$:  
the latter rely chiefly on properties of conformal mappings, while the 
former are based entirely on particular forms of the comparison 
principle.  
Our discussion above for a trihedral angle applies only to the case 
$H \ge 0$.  But Theorem~\ref{t1} shows that under the hypotheses 
of that theorem there can be no drop with $H \ge 0$ in a wedge.  Consider 
as an example a trihedral angle formed by three orthogonal planes, 
with the same contact angle $\gamma = \cos^{-1}(\sqrt{3}/3)$ on all planes. 
A symmetrically placed planar surface ${\cal S}$ then determines a 
drop that wets all planes, and by Theorem~\ref{t3} it is 
the unique such drop for which $H = 0$. If we now rotate one of the 
planar sides $\Pi$ into the angle about its intersection line with 
one of the other planes, then a surface ${\cal S}$ with $H > 0$ can be 
found as a spherical cap, and according to Theorem~\ref{t3} it is the 
unique disk-type surface with the given $\gamma$ and that value of $H$.  
On the other hand, if we rotate $\Pi$ an amount less than $\pi/2$ in 
the other direction then $H < 0$ for the spherical cap, and no 
uniqueness theorem is available.  But if we continue the rotation 
back until $\Pi$ coincides with the other plane and the trihedral 
angle becomes a wedge, then Theorem~\ref{t1} guarantees the uniqueness 
of the spherical cap among all competing disk-type surfaces.  On the 
basis of these remarks, we formulate
\begin{conjecture}
Theorem~\ref{t3} holds without the hypothesis $H \ge 0$.
\end{conjecture}

\appendix
\section*{Appendix}
The method of this paper (for the case $V = 2$) yields as corollary 
a conceptually simpler proof of a theorem of Nitsche \cite{NitSta},  
that an immersed disk ${\cal S}$ of constant mean curvature that 
intersects a sphere $\Sigma$ at constant angle along a smooth closed 
curve ${\cal C}$ is necessarily a spherical cap.  In fact, every curve is a 
curvature line on $\Sigma$, thus ${\cal C}$ has that property and thus by 
\cite{JoaDem}  
${\cal C}$ is a curvature line on ${\cal S}$.  We map ${\cal S}$ conformally 
to a unit disk in the plane, and the disk to a strip by (\ref{2}).  We 
find immediately using (\ref{6}) that the second fundamental form 
vanishes identically in the strip, from which follows that 
${\cal S}$ is totally umbilic and hence a metric sphere. 

\bigskip

\bigskip

{\small We wish to thank Erich Miersemann for a number of helpful 
discussions, and notably for informing us of the paper of 
Azzam \cite{AzzBeh}.}  

\bibliographystyle{alpha}  
\bibliography{bib/drops}

\end{document}